\def\hook{ \, \mbox{\raisebox{-.19ex}{\rule[.1mm]{2.15mm}
{.1mm}}\raisebox{-.15ex}{\rule{0.1mm}{2.3mm}}}~}
\newtheorem{main}{Theorem}

\newtheorem{thm}{Theorem}
\newtheorem{prop}{Proposition}
\newtheorem{lem}{Lemma}
\newtheorem{cor}{Corollary}
\newtheorem{defn}{Definition}
\newenvironment{proof}{\medskip
\noindent {\bf Proof.}}{\hfill \rule{.5em}{1em} \\}
\newenvironment{example}{\medskip
\noindent {\bf Example.}}{}

\documentstyle[12pt]{article}
\let\Bbb\bf
\title{A Finiteness Theorem for Quaternionic-K\"ahler Manifolds with
Positive Scalar Curvature}
\author{  Claude LeBrun\thanks{Supported
in part by  NSF grant DMS-9003263.
}\\SUNY Stony
 Brook}
\date{}
\begin{document}
\maketitle
\abstract{
We study the topology and geometry
of those   compact  Riemannian manifolds $(M^{4n},g)$, $n \geq 2$,
with positive scalar curvature and  holonomy
in $Sp(n)\cdot Sp(1)$. Up to homothety,
 we show that
there are only finitely many such manifolds of any  dimension $4n$.
}
\vfill
\pagebreak
\hfill
\begin{minipage}[t]{3.5in}
\begin{verse} \em
\ldots che tu mi sie di tuoi prieghi cortese,\\
in Fano, s\`{\i} che ben per me s'adori\\
pur ch'i' possa purgar  le gravi offese.\\
\end{verse} \hfill
\begin{minipage}[t]{2.5in}
\begin{quote}
  Dante Alighieri\\ \em Purgatorio V, 70-72
\end{quote}
\end{minipage}
\end{minipage}
\bigskip

\bigskip

Let $(M^{\ell},g)$ be a connected Riemannian manifold. If $x\in M$ is an
arbitrary base-point, one defines the {\em holonomy group}
\linebreak ${\cal H}(M,g,x)
\subset {\rm End}(T_xM)$ of $(M,g,x)$
to be the set of linear maps $T_xM\to T_xM$ obtained by Riemannian
parallel transport of tangent vectors around piece-wise smooth loops
based at $x$; the {\em reduced holonomy group} ${\cal H}_0(M,g,x)$ is similarly
defined, but now
using only null-homotopic loops. The latter is automatically a
connected Lie subgroup
of the orthogonal transformations of the tangent space $T_xM$,
and so may be identified with
a Lie subgroup of $SO(\ell)$ by choosing an orthonormal basis for $T_xM$;
and the conjugacy
class of this subgroup is uneffected by a
change of basis or  base-point $x$.
 As was first pointed out by Berger,
relatively few groups can actually arise in this way.
If we exclude the so-called {\em locally reducible} manifolds, meaning those
which  are locally Riemannian Cartesian products of lower-dimensional
manifolds, and the {\em locally symmetric} manifolds, meaning those
for which the curvature tensor is covariantly constant, the only
possibilities \cite{besse}\cite{Sbook}
for the reduced holonomy of $M^{\ell}$ are those which
appear on the following list:
\bigskip

\begin{center}
\begin{tabular}{|c||c|c|}  \hline
$\ell$ & ${\cal H}_0$ & geometry
\\ \hline \hline
$\ell$ & $SO(\ell)$&generic\\ \hline
$2m\geq 4$ &$U(m)$ & K\"ahler \\ \hline
$2m \geq 4$&$SU(m)$& Ricci-flat K\"ahler\\ \hline
$4n\geq 8$ & $Sp(n)$& hyper-K\"ahler\\ \hline
$4n\geq 8$ &$Sp(n)\times Sp(1)/{\Bbb Z}_2$& quaternionic-K\"ahler\\ \hline
7 &$G_2$& imaginary Cayley\\ \hline
8&$Spin(7)$& octonionic\\ \hline
\end{tabular}
\end{center}
In particular, the universal cover of any complete Riemannian manifold is
the  product of  globally symmetric spaces and manifolds whose
holonomy appears on this list. Each of these possible holonomies
may therefore be thought of as representing a ``fundamental geometry,''
each of which has its own peculiar flavor.

Of these fundamental geometries, the  {\em quaternionic-K\"ahler}
or  \linebreak $Sp(n)\cdot Sp(1):=
Sp(n)\times Sp(1)/{\Bbb Z}_2$ possibility in in some ways the most
enigmatic. It is  the only holonomy geometry which
forces the metric to be Einstein but not Ricci-flat. In particular,
the scalar curvature of  a quaternionic-K\"ahler manifold is constant, and the
sign of this constant turns out to have a decisive influence
on the nature  of manifold in question.
Thus, while complete, non-compact, non-symmetric  quaternionic-K\"ahler
manifolds of negative scalar curvature
 exist in profusion
\cite{a}\cite{g}\cite{L2}, a complete  quaternionic-K\"ahler
manifold of positive scalar curvature must
 be compact, and the only  known examples
 of such manifolds are symmetric.
Indeed,  Poon and Salamon \cite{PS},
 generalizing earlier work of Hitchin \cite{H},
 have proved that
there are no others in dimension 8, and it is therefore
tempting to conjecture that this situation persists in higher dimensions.
This article  will
indicate some  new evidence supporting such a
conjecture. In particular,  a complete proof  the following
result is given:

\begin{main}[Finiteness Theorem]
For any $n$, there are, modulo isometries and rescalings,
 only finitely many compact quaternionic-K\"ahler $4n$-manifolds of
positive scalar
curvature.
\end{main}

Of course, this doesn't predict that every such manifold is symmetric.
But the other main result which will be described
in detail herein does allow one to draw such a conclusion
 for certain
topological types:

\begin{main}[Strong Rigidity]
 Let $(M,g)$ be a compact quaternionic-K\"ahler manifold of
positive scalar curvature. Then $\pi_1(M)=0$ and
$$\pi_2(M)= \left\{
\begin{array}{cl}
0&(M,g)={\Bbb HP}_n\\
{\Bbb Z}&(M,g)=Gr_2({\Bbb C}^{n+2})\\
\mbox{finite}\supset {\Bbb Z}_2&\mbox{otherwise.}
\end{array}
\right. $$
\end{main}

It might be particularly emphasized that these result rely very heavily on
the global aspects of the problem. In particular, if one chooses to consider
compact {\em orbifolds} rather than manifolds, both fail
\cite{gl}, even though the more general metrics involved are still
 smooth everywhere---
at least when viewed from  the skewed perspective of
certain multi-valued coordinate charts.

Many of the results described herein were obtained as part of
a joint research project with Simon Salamon, and further details
will appear in our joint paper \cite{ls}.

\section{Preliminaries}
\begin{defn} Let $(M, g)$ be a connected Riemannian $4n$-manifold, $n\geq 2$.
We will say that $(M,g)$ is a quaternionic-K\"ahler manifold iff
the holonomy group ${\cal H}(M,g)$ is conjugate to $H\cdot Sp(1)$
for some Lie subgroup $H\subset Sp(n)\subset SO (4n)$.
\end{defn}

\begin{example} The quaternionic projective spaces
 $${\Bbb HP}_n=Sp(n+1)/(Sp(n)\times
Sp(1))$$ are
quaternionic-K\"ahler manifolds. So are the complex Grassmannians
$$Gr_2({\Bbb C}^{n+2})=SU(n+2)/S(U(n)\times U(2))$$ and the
oriented real Grassmannians $$\tilde{Gr}_4({\Bbb R}^{n+4}) =
SO(n+4)/(SO(n)\times SO(4)).$$
\end{example}

In fact, these examples very nearly exhaust the compact homogeneous
examples of quaternionic-K\"ahler manifolds. Indeed \cite{a2},
every  such homogeneous space is a symmetric space, and \cite{wo} there is
exactly one such symmetric space for each compact simple Lie algebra.
They can be constructed as follows: let $G$ be a compact simple
centerless group, and let $Sp(1)$ be mapped to $G$ so that
its root vector is mapped to a root of highest weight. If
$H$ is the centralizer of this $Sp(1)$, then the symmetric space
$M=G/(H\cdot Sp(1))$ is quaternionic-K\"ahler, and every
compact homogeneous quaternionic-K\"ahler manifold arises this way.

How typical are these symmetric examples? One geometric feature of
any irreducible symmetric space is that it must be Einstein,
with non-zero scalar curvature.   This, it turns out, also happens
for quaternionic-K\"ahler manifolds:

\begin{prop}[Berger] Every quaternionic-K\"ahler manifold is
Einstein, with non-zero scalar curvature. \end{prop}

\noindent For details, see \cite{besse}. In particular, a complete
quaternionic-K\"ahler manifold has constant scalar curvature.

\begin{defn} We will say that a quaternionic-K\"ahler manifold is
{\em positive} if it is complete and has positive scalar curvature.
\end{defn}

\noindent
It is now an an immediate consequence of  Myers' theorem that
a positive quaternionic-K\"ahler manifold is compact
and has finite fundamental group.
Unfortunately, however, the only known
 positive quaternionic-K\"ahler manifolds
are the previously mentioned symmetric spaces! The main objective of the
present
article will be to explain why this situation is hardly surprising.
The main tool in our investigation will be the next result:

\begin{thm}[Salamon {\rm \cite{S}}; B\'erard-Bergery {\rm \cite{beber}}]
Let $(M^{4n},g)$ be a quaternionic-K\"ahler manifold. Then there is
a complex manifold $(Z,J)$ of complex dimension $2n+1$,
called the {\em twistor space} of $(M,g)$,  such that
\begin{itemize}
\item there is a smooth fibration $\wp : Z\to M$ with fiber $S^2$;
\item each fiber of $\wp$ is a complex curve in $(Z,J)$ with
normal bundle holomorphically isomorphic to $[{\cal O}(1)]^{\oplus 2n}$,
where ${\cal O}(1)$ is the point-divisor line bundle on ${\Bbb CP}_1$; and
\item there is a complex-codimension 1 holomorphic sub-bundle $D\subset TZ$
which is maximally non-integrable and transverse to the fibers of $\wp$.
\end{itemize}
Moreover, if $(M,g)$ is {\em positive},
then $Z$ carries a K\"ahler-Einstein metric of
positive scalar curvature such that
\begin{itemize}
\item $\wp$ is a Riemannian submersion;
\item $D$ is the orthogonal complement of the vertical tangent bundle of $\wp$;
and
\item the induced  metric on each fiber of $\wp$  has constant curvature.
\end{itemize}
If $(M,g)$ is instead negative, there is an {\em indefinite}
K\"ahler-Einstein pseudo-metric on $Z$ with all these properties.\label{sal}
\end{thm}

In particular, the twistor space $Z$ of a  positive
quaternionic-K\"ahler manifold is  {\em Fano}:

\begin{defn} A Fano manifold is a compact complex manifold
$Z$ such that $c_1(Z)$ can be represented by a positive (1,1)-form.
\end{defn}

\noindent
That is, a Fano manifold is a compact complex manifold which admits
K\"ahler metrics of positive Ricci curvature.
Every Fano manifold is simply connected, since $c_1>0~\Rightarrow ~\chi ({\cal
O})
=h^0 ({\cal O})
=1$ by the Kodaira vanishing theorem, thus forbidding the possibility that
the manifold might have a finite cover. Applying the exact homotopy
sequence of $Z\to M$, we now  conclude the
 following:

\begin{prop} Any positive quaternionic-K\"ahler manifold is compact
and simply connected.
\end{prop}

A completely different and extremely important feature of our
twistor spaces is the holomorphic hyperplane distribution $D$,
which gives a so-called {\em complex contact structure} to $Z$.
Such structures will be  discussed systematically in \S \ref{cont}.

Our definition of quaternionic has carefully avoided the case of $n=1$;
after all, $Sp(1)\cdot Sp(1)$ is all of $SO(4)$, so such a  holonomy
restriction says nothing at all. Instead, we choose the our definition
in order to insure that Theorem \ref{sal} remains valid:

\begin{defn} A Riemannian 4-manifold $(M,g)$ is called
 quaternionic-K\"ahler if it is Einstein, with non-zero scalar curvature,
and half-conformally flat.
\end{defn}

\section{Complex Contact Manifolds} \label{cont}

\begin{defn} A complex contact manifold is a pair $(X,D)$, where
$X$ is a complex manifold  and
$D\subset TX=T^{1,0}X$ is a codimension-one holomorphic sub-bundle
which is maximally non-integrable in the sense that the
 O'Neill tensor
\begin{eqnarray*}
D\times D &\to& TX/D\\
(v,w)&\mapsto & [v,w]\bmod D
\end{eqnarray*}
is everywhere  non-degenerate.
\end{defn}

\begin{example} Let $Y_{n+1}$ be any complex manifold, and let
$X_{2n+1}={\Bbb P}(T^{\ast}Y)$
be its projectived
 holomorphic cotangent bundle;  dually stated,  $X$ is the  Grassmann bundle
of complex $n$-planes in $TY$. Let $\pi :X\to Y$ be the canonical projection,
and let $D\subset TX$ be the sub-bundle defined by $D|_P:=\pi^{-1}_{\ast}(P)$
for all complex n-planes $P\subset TY$. Then $D$ is a complex contact structure
on $X$.
\end{example}\bigskip

The condition of non-integrability has a very useful reformulation,
which we shall now describe. Given a codimension-one holomorphic sub-bundle
$D\subset TX$, let $L:=TX/D$ denote the quotient line bundle.
Letting $\theta : TX\to L$ be the
tautological  projection, we may think of $\theta$ as a line-bundle-valued
1-form
$$\theta \in \Gamma (X, \Omega^1(L))~,$$
and so attempt to form its exterior derivative $d\theta$.
Unfortunately, this ostensibly depends on a choice of local
trivialization; for if $\vartheta$ is any 1-form,
$d(f\vartheta)=fd\vartheta +df\wedge\vartheta$. However, it is now clear that
$d\theta|_D$ {\em is} well defined as a section of $L\otimes \wedge^2D^{\ast}$,
and an elementary computation, which we leave to the reader, shows that
 $d\theta|_D$, thought of in this way, is exactly the  O'Neill tensor mentioned
above. Now if the skew form $d\theta|_D$ is to be non-degenerate, $D$ must
have positive even rank $2n$, so that  $X$ must have odd complex dimension
$2n+1\geq 3$.
Moreover, the non-degeneracy exactly requires that
$$\theta \wedge (d\theta )^{\wedge n}\in \Gamma (X, \Omega^{2n+1}(L^{n+1}))$$
is nowhere zero. But this provides a bundle isomorphism between
$L^{\otimes (n+1)}$
and the anti-canonical line bundle $\kappa^{-1}=\wedge^{2n+1}T^{1,0}X$.

Conversely, let $X$ be a simply-connected  compact complex $(2n+1)$-manifold,
and suppose that $c_1(X)$ is divisible by $n+1$. Then there is a
unique holomorphic line bundle $L:=\kappa^{-1/(n+1)}$ such that
$L^{\otimes (n+1)}\cong \kappa^{-1}$. If we are then given a
twisted holomorphic 1-form
$$\theta \in \Gamma (X, \Omega^1(\kappa^{-1/(n+1)}))$$
we may then construct
$$\theta \wedge (d\theta )^{\wedge n}\in \Gamma (X,
\Omega^{2n+1}(\kappa^{-1}))=
\Gamma (X,{\cal O})= {\Bbb C}~.$$
If this constant is non-zero,  $D=\ker \theta$ is then a complex
contact structure.

This simple observation has powerful consequences:

\begin{prop} Let $X_{2n+1}$ be a simply connected compact complex
manifold, and let $\cal G$ denote  the identity component of the
group of biholomorphisms $X\to X$. Then $\cal G$ acts transitively
on the set of complex contact structures on $X$. \label{trans}
\end{prop}

\begin{proof} We may assume that there is  at least one
complex contact structure on $X$, since otherwise there is nothing to prove.
In this case, the canonical line bundle $\kappa$ has a root $\kappa^{1/(n+1)}$,
and there is only one such  root  because $H^1(X, {\Bbb Z}_{n+1})=0.$
Thus any complex contact structure is determined by a
class $[\theta]\in
{\Bbb P }\Gamma (X, \Omega^1\otimes (\kappa^{-1/(n+1)}))$
satisfying $\theta\wedge (d\theta)^n\neq 0$.
The group $\cal G$ acts on this
projective space $\cong {\Bbb P}_m$ in a manner preserving the
hypersurface $S$ defined by  $\theta\wedge (d\theta)^n= 0$,
and so  partitions ${\Bbb P}_m-S$ into orbits; since ${\Bbb P}_m-S$
is connected,
it therefore suffices
to prove that each orbit is open, and for this it would be enough to
prove that  the holomorphic
vector fields generating the action of
the the Lie algebra of $\cal G$ on ${\Bbb P}_m-S$
span the tangent space at each point.

To prove the last statement, let $\theta\in \Gamma
(X, \Omega^1\otimes (\kappa^{-1/(n+1)}))$ be any contact form,
and let $\phi \in \Gamma
(X, \Omega^1\otimes (\kappa^{-1/(n+1)}))$ be any other section.
If $D$ denotes the kernel of $\theta$, $d\theta |_D: D\to D\otimes
 \kappa^{-1/(n+1)}$ is an isomorphism of holomorphic vector bundles,
so we can define a holomorphic vector field $v\in \Gamma
(X, {\cal O} (D))$ by $v=(d\theta |_D)^{-1}(\phi)$. We then
have $\pounds_{\xi} \theta\equiv {\xi} \hook d\theta\equiv \phi \bmod\theta$,
so that action of the Lie algebra of $\cal G$ spans the
tangent space of ${\Bbb P }\Gamma (X, \Omega^1\otimes (\kappa^{-1/(n+1)}))$
at $[\theta ]$, thus proving the proposition. \end{proof}

\begin{cor}
Two simply-connected compact complex manifolds are
complex-contact isomorphic iff the underlying
complex manifolds are biholomorphically equivalent.
\end{cor}

This will now yield a  result which is crucial for our purposes.

\begin{defn} We will say that
two Riemannian
manifolds $(M_1,g_1)$ and $(M_2,g_2)$ are
{\em homothetic} if there exists a diffeomorphism
$\Phi : M_1\to M_2$ such that $\Phi^{\ast}g_2=cg_1$ for some
 constant $c>0$. Such a map $\Phi$ will be called a
{\em homothety}.
\end{defn}

\begin{prop} Two
positive quaternionic-K\"ahler manifolds
are
homothetic iff their twistor spaces
are biholomorphic. \label{perfect}
\end{prop}

\begin{proof}
\label{final} Let $(M, g)$ and $(\tilde{M}, \tilde{g})$ be two
given quaternionic-K\"ahler manifolds,  $\wp :Z\to M$ and
$\tilde{\wp} :\tilde{Z}\to\tilde{M}$  their twistor spaces,
$h$ and $\tilde{h}$ the K\"ahler-Einstein metrics of
$Z$ and $\tilde{Z}$. We also
 suppose that a biholomorphism $\Phi: Z\to \tilde{Z}$
is given to us.
For some positive constant $c>0$,
$h$ and $c\tilde{h}$ have the same scalar curvature;
and notice that replacing $\tilde{h}$ with  $c\tilde{h}$ just corresponds to
replacing $\tilde{g}$ with $c\tilde{g}$.
Now $\Phi^{\ast}c\tilde{h}$ is a K\"ahler-Einstein metric
on $Z$ with the same scalar curvature as $h$, and
 the Bando-Mabuchi theorem \cite{BM} on the
uniqueness of K\"ahler-Einstein metrics
now asserts that there exists a biholomorphism
$\Psi: Z\to Z$ such that $\Psi^{\ast}(\Phi^{\ast}c\tilde{h})=h$.

Let $N\subset\Gamma (Z, \Omega^1(\kappa^{-1/(n+1)}))$
be  defined by $\phi\wedge (d\phi)^{\wedge n}=1$.
 Proposition \ref{trans} implies that
 a finite connected cover $\cal G$ of the connected component
of the automorphism group of $(Z, J)$ acts transitively on
$N$, since, in the notation of the proof of that
proposition, $N\to {\Bbb P}_m-S$ is a finite covering.
Because $h$ is K\"ahler-Einstein, with
positive scalar curvature, the Killing fields are a real form of the
algebra of holomorphic vector fields, and a finite cover $G$ of
the  connected component of the isometry group  of $(Z,h)$ is therefore
a compact real form of $\cal G$.  Morse theory
now predicts  that one orbit
of the action of $G$ on $N$ is precisely the set of
critical points of the $G$-invariant strictly plurisubharmonic
function $f:N\to {\Bbb R}$ given by $\phi\mapsto \|\phi\|^2_{L^2,h}$.
On the other hand, the derivative of $f$ at a contact form $\phi$
in the direction of a real-holomorphic vector field $\xi$ on $Z$
is given by
$$df(\xi )|_{\phi}={\textstyle \frac{1}{(2n)!}}
\int_X d (\xi \hook \omega)
\wedge |\phi|^2_h\left(\beta_{\phi}-{\textstyle \frac{n+2}{n+1}}
\omega^{2n}\right)~,$$
where $\omega$ is the K\"ahler form of $(Z,J,h)$ and
$\beta_{\phi}$ is the $(2n,2n)$ form obtained by
orthogonally extending the restriction of
$\omega^{2n}$ from $D=\ker \phi$  to $TZ$. For
 the canonical contact form $\theta$ associated
with the quaternionic-K\"ahler metric $g$ by the twistor construction,
$|\theta |_h$ is constant and $\beta_{\theta}=\wp^{\ast}(2n)!
d\mbox{\rm vol}_g$  is closed, so that
$\theta$ is a
critical point of $f$; but the    same argument    applies  equally to
the contact from of $\tilde{Z}$, and hence to
the pull-back of this contact
form  via the holomorphic isometry $\Phi\Psi$.
Hence there is a  holomorphic
isometry $\Xi: Z\to Z$ sending the first of these contact structures to the
second, and $\Phi\Psi\Xi : Z\to \tilde{Z}$ is a then biholomorphism  which
sends
$h$ to $c\tilde{h}$ and $D$ to $\tilde{D}$. Since the vertical
tangent spaces  of
$\wp$ and $\tilde{\wp}$ are the orthogonal complements of
$D$ and $\tilde{D}$ with respect to $h$ and $\tilde{h}$, respectively,
it  follows that $\Phi\Psi\Xi$ sends fibers of
$\wp$ to fibers of $\tilde{\wp}$, and so covers a diffeomorphism
$F:M\to \tilde{M}$. Moreover, since the $\wp$ and $\tilde{\wp}$
are Riemannian submersions, one has $F^{\ast}c\tilde{g}=g$,
and $F$ is thus a homothety between $(M,g)$ and $(\tilde{M}, \tilde{g})$.
\end{proof}


For less precise but more broadly applicable theorems on the invertibility
of the twistor construction, cf. \cite{L0}\cite{BE}.

\begin{defn} Let $(X_{2n+1},D)$ be a complex contact
manifold. An $n$-dimensional submanifold $\Sigma_n\subset X_{2n+1}$ is called
{\em Legendrian} if $T\Sigma\subset D$.
\end{defn}

\begin{lem} Let $(X_{2n+1},D)$ be a complex contact
manifold, and let \linebreak
$\pi: X\to Y_{n+1}$ be a proper holomorphic submersion with
Legendrian fibers. Then $X\cong {\Bbb P}(T^{\ast} Y)$ as complex
contact manifolds. \label{class}
\end{lem}
\begin{proof} Define $\Psi : X\to Gr_{n}(TY)$ by $x\to \pi_{\ast}(D_x)$.
This map preserves the contact structure; and
since the pull-back   of the contact form of $Gr_{n}(TY)$ via $\Psi$
is the contact form of $X$, $\Psi^{\ast}$ induces an isomorphism
between forms of top degree. In other words, $\Psi$ is a
submersion onto its
image,
and, in particular, induces a submersion from each fiber of
$X$ onto its image in the  fiber of $Gr_{n}(TY)={\Bbb P}(T^{\ast} Y)$.
By the  properness assumption, $\Psi$ is  fiber-wise
therefore a covering map.
But the fibers of ${\Bbb P}(T^{\ast} Y)$ are projective spaces,
and so simply connected. Hence $\Psi$ is an injective holomorphic
submersion, and so
biholomorphic.
\end{proof}

\begin{defn} If $(X,D)$ is a complex contact manifold such that $X$ is Fano,
we will say that $(X,D)$ is a
Fano contact manifold.
\end{defn}

\begin{lem} Let $\varpi: {\cal X} \to {\cal B}$ be
a  holomorphic family of Fano contact manifolds
with smooth connected parameter space--- that is, let $\cal B$
be a connected complex manifold,  $\varpi$ a proper holomorphic
submersion with Fano fibers, and
assume that $\cal X$ is equipped with a maximally
non-integrable, complex codimension 1
sub-bundle $D\subset \ker \varpi_{\ast}$ of the vertical tangent bundle.
 Then any two fibers
$(X_0, D|_{X_0})$  and $(X_t,D|_{X_t})$
are isomorphic as complex contact manifolds.\label{rig}
\end{lem}
\begin{proof}
 Since any two points in $\cal B$ can be joined by
a finite chain of holomorphic images of the unit disk $\Delta \subset {\Bbb
C}$,
it suffices to prove the lemma when the base $\cal B$
is a disk $\Delta$.

We now proceed as in \cite{L1}.
By Darboux's theorem, any complex contact structure
in dimension $2n+1$ is locally complex-contact isomorphic to the one on
 ${\Bbb C}^{2n+1}$ determined by the 1-form
$$\vartheta = dz^{2n+1}+\sum_{j=1}^{n}z^jdz^{n+j}~, $$
so we may cover our family $\varpi :{\cal X}\to \Delta$ by Stein sets
$U_j$ on which we have
holomorphic charts $\Phi_j: U_j\hookrightarrow {\Bbb C}^{2n+1}\times \Delta$
such that the last coordinate is given by $\varpi$ and the fiber-wise contact
structure on $\cal X$ agrees with that induced by $\Phi_j^{\ast} \vartheta$.
Letting $t$ denote the standard complex coordinate on
$\Delta$, we lift $d/dt$ to each $U_j$ as the vector field
$v_j(\Phi_j^{-1})_{\ast}d/dt$, and observe that the $t$-dependent
vertical vector field $w_{jk}:=v_j-v_k$ satisfies $\pounds_{w_{jk}}\theta
\propto \theta$.

Let $f_{jk}:=\theta (w_{jk})\in \Gamma (U_j\cap U_k, {\cal O} (L))$,
and notice that the
 collection $\{ f_{jk}\}$ is  a \v{C}ech cocycle representing an element
of   $H^1({\cal X}, {\cal O} (L))$. On the other hand, since $L$ is
a fiber-wise $(n+1)^{st}$-root of the vertical anti-canonical
bundle $\kappa^{-1}$,
and since each fiber $X_t$ of $\varpi$ is assumed to be a Fano manifold,
the bundle $\kappa^{-1}\otimes L$ is fiber-wise positive, and
 $H^1(X_t, {\cal O}(L))=0$ $\forall t\in \Delta$ by the
Kodaira vanishing theorem. Thus the first direct image sheaf
$\varpi_{\ast}^1{\cal O}(L)$ is zero. Since $\Delta$ is Stein,
the Leray spectral sequence now yields  $H^1({\cal X}, {\cal O}(L))=0$.
Hence there exist sections $h_j\in  \Gamma (U_j\, {\cal O} (L))$
such that $f_{jk}=h_j-h_k$ on $U_j\cap U_k$.

On $U_j$ there is now a unique vertical holomorphic vector field $u_j$ such
that
$\theta (u_j)=h_j$ and $\pounds_{u_j}\theta\propto\theta$. Indeed,
taking a local trivialization of $L$ so as to locally represent
$\theta$ by a holomorphic 1-form $\vartheta$, a vector field
$u$ satisfies $\pounds_u\theta\propto\theta$ iff
$$u\hook d\vartheta \equiv -d(u\hook \vartheta ) \bmod \vartheta~,$$
so that such a field is uniquely determined by an arbitrary local function
$f=\vartheta (u)= u\hook \vartheta$. We therefore conclude that
$v_j-v_k=w_{jk}=u_j-u_k$ on $U_j\cap U_k$, and the vector field
$v=v_j-u_j$ is therefore globally defined.
Since $\pounds_{v_j}\vartheta =0$ and $\pounds_{u_j}\vartheta \equiv 0
\pmod{\vartheta, dt},$  the flow of $v=v_j-u_j$ preseves the fiberwise contact
structure on $\cal X$.
And since   $\varpi$ is a proper map, we can now integrate the flow
of our lift $v$ of $d/dt$ to produce a fiber-wise contact
biholomorphism between $\cal X$ and $ X_0\times \Delta$.
In particular, any fiber $X_t$ is complex-contact equivalent to the central
fiber
$X_0$.
\end{proof}

\section{Mori Theory}\label{this}

Mori's theory of extremal rays \cite{M} has
led to a startling series of advances in the classification of
 complex algebraic varieties, especially in the Fano case which
interests us. One beautiful consequence of this is the
so-called {\em contraction theorem}: if $X$ is a Fano manifold,
there is always a map $\Upsilon : X\to Y$ to some other  variety $Y$
 which decreases
the second Betti number $b_2$ by one, and where the kernel
of $ \Upsilon_{\ast} :H_2(X, {\Bbb R})\to H_2(Y, {\Bbb R})$
 is generated by the class of a rational
holomorphic curve ${\Bbb CP}_1\subset X$. (The positive half
of such a one-dimensional subspace $\ker \Upsilon_{\ast}\subset H_2(X, {\Bbb
R})$
is called an ``extremal ray''.)
If $b_2(X)=1$, this tells us next to nothing,
because we can take $Y$ to be a point; but for $b_2(X)\geq 2$, it is
quite a powerful tool. In particular, it gives rise to the following
very useful result of Wi\'sniewski \cite{W}:

\begin{thm}[Wi\'sniewski]
Let $X$ be a Fano manifold of dimension
$2r-1$ for which $r|c_1$. Then $b_2(X)=1$ unless $X$ is one of the
following:
(i) ${\Bbb CP}_{r-1}\times Q_r$\ ; (ii) $ {\Bbb P}(T^{\ast}{\Bbb CP}_r)$\ ; or
(iii) ${\Bbb CP}_{2r-1}$ blown up along ${\Bbb CP}_{r-2}$.\label{wis}
\end{thm}

Here $Q_r\subset {\Bbb CP}_{r+1}$ denotes the r-quadric,
 while the projectivization
of a bundle $E\to Y$ is defined by ${\Bbb P}(E):=(E-0_Y)/({\Bbb C}-0)$.
The essence  of the proof is that, since the
rational curves collapsed by the Mori contraction
 have, in these circumstances,  normal bundles of rather large index,
they are so mobile that they sweep out projective spaces
of comparatively large dimension, and these must therefore
be the fibers of the contraction map.

The following is now an easy  consequence:

\begin{cor} Let $(X_{2n+1},D)$ be a Fano contact manifold. If \linebreak
$b_2(X)>1$, then $X= {\Bbb P}(T^{\ast}{\Bbb CP}_{n+1})$. \label{boop}
\end{cor}
\begin{proof}
Setting $r=n+1$, we notice that the existence of a
contact structure implies that $(n+1)|c_1$. We may
therefore invoke  Theorem \ref{wis}.
On the other hand, spaces (i) and (iii) aren't complex contact manifolds,
since
$\Gamma ({\Bbb CP}_{r-1} , \Omega^1 (1))=0$ and therefore the
obvious foliations by
${\Bbb CP}_{r-1}$'s would necessarily have Legendrian leaves, implying
(by Lemma \ref{class}) that
these spaces would then have to be of the form ${\Bbb P}(T^{\ast}Y)$,
where
$Y$ is the leaf space  $Q_r$ or ${\Bbb CP}_{r}$--- a
contradiction. So the only candidate left is (ii), and this {\em is} in fact
a contact manifold.
\end{proof}

\begin{thm}{\rm \cite{L3}}
Let $(M,g)$ be a compact
quaternionic-K\"ahler $4n$-manifold with $s>0$.
Then either
\begin{description}
\item{(a)} $b_2(M)=0$; or else
\item{(b)}  $M=Gr_2({\Bbb C}^{n+2})$
with its symmetric-space metric.\label{grass}
\end{description}
\end{thm}
\begin{proof}
By  the Leray-Hirsch theorem on sphere bundles, the second Betti numbers of
$M^{4n}$ and its twistor
space  $Z_{2n+1}$ are related by  $b_2(Z)=b_2(M)+1$.
Since $Z$ is a Fano contact manifold,
 $b_2(M)>0\Rightarrow Z={\Bbb P}(T^{\ast}{\Bbb CP}_{n+1})$ by
Corollary \ref{boop}.
But this is the twistor space of $Gr_2({\Bbb C}^{n+2})$. The result
therefore follows by Proposition \ref{perfect}.
\end{proof}

\setcounter{main}{1}
\begin{main}[Strong Rigidity]
 Let $M$ be a compact quaternionic-K\"ahler manifold of
positive scalar curvature. Then $\pi_1(M)=0$ and
$$H_2(M, {\Bbb Z} )= \left\{
\begin{array}{cl}
0&M={\Bbb HP}_n\\
{\Bbb Z}&M=Gr_2({\Bbb C}^{n+2})\\
\mbox{finite}\supset {\Bbb Z}_2&\mbox{otherwise.}
\end{array}
\right. $$\label{tor}
\end{main}
\begin{proof} If $(M, g)$ is not homothetic to  the symmetric space
$Gr_2({\Bbb C}^{n+2})=SU(n+2)/S(U(n)\times U(2))$,
$b_2(M)=0$
by Theorem \ref{grass}, so that
$H^2(M, {\Bbb Z})=0$ and
$H_2(M, {\Bbb Z})$ is  finite. Since we also know that $H_1(M, {\Bbb Z})=0$,
$H^2(M, {\Bbb Z}_2)$ is exactly the the 2-torsion  of
$H_2(M, {\Bbb Z})$ by the universal coefficients theorem.
If, on
 the other hand, $(M, g)$ is not homothetic to
the symmetric space  ${\Bbb HP}_n$, the class $\varepsilon\in H^2(M, {\Bbb
Z}_2)$
must be non-zero \cite{S}, and the finite group
$\pi_2(M)=H_2(M, {\Bbb Z})$ must therefore contain an element of order 2.
\end{proof}

\section{The Finiteness Theorem}

\begin{thm} Up to biholomorphism, there are only finitely many
Fano contact manifolds of any given dimension $2n+1$.\label{fanite}
\end{thm}
\begin{proof}
By Wisniewski's theorem, we may restrict our attention
to Fano manifolds with  $b_2=1$. A theorem of Nadel \cite{nad}
then asserts that there are only a finite number of deformation
types of any fixed dimension.\footnote{It is now known \cite{kmm} that this
is true even {\em without} the  restriction $b_2=1$.}

For any fixed deformation type, we may embed each Fano manifold
in a fixed projective space ${\Bbb CP}_N$ in such a manner that
the restriction of the generator $\alpha \in H^2({\Bbb CP}_N, {\Bbb Z})$
is a fixed multiple $\ell c_1(Z)/q$ of
the anti-canonical class, and
 we may freely choose the positive integers  $\ell$ and $q$ as
long as
$q|c_1$  and  $\ell$ is sufficiently large.
 Thus, let $\cal F$ denote the set of
all complex submanifolds  $Z\subset {\Bbb CP}_N$ of some fixed
dimension $m$ and degree $d$,
and with the additional property that, for  fixed integers
$\ell, q$,  the restriction of the hyperplane class
 is $\ell c_1(Z)/q$. Thus
$\cal F$ is a Zariski-open subset in a component of the
Chow variety, and so, in particular, is quasi-projective.
There is now a tautological family
\begin{eqnarray*}
{\cal Z} & \hookrightarrow & {\cal F}\times {\Bbb CP}_N    \\
{\varpi}\downarrow& & ~~~~\downarrow  \\
{\cal F}   & ~~=  & ~~~~{\cal F}
\end{eqnarray*}
such that the fiber of $Z\in {\cal F}$ is the submanifold
$Z\subset {\Bbb CP}_N$.

We now assume moreover that the dimension $m$  is an odd number $2n+1$,
take $q=n+1$, and choose $\ell \gg 0$ such that $\gcd (n+1, \ell)=1$.
Letting $V\to {\cal Z}$ denote the vertical tangent bundle
$\ker \varpi_{\ast}$, the  vertical anti-canonical line bundle
$\kappa^{-1}:= \wedge^{2n+1}V$
has a consistent fiber-wise
$(n+1)^{st}$-root $L\to {\cal Z}$; indeed, using the Euclidean algorithm to
write
$1=a(n+1)+b\ell$, we may define $L$ by $L=\kappa^{-a}\otimes {\cal H}^{b}$,
where $\cal H$ is the pull-back of ${\cal O}(1)$ from ${\Bbb CP}_N$
to $\cal Z$.

Let ${\cal F}_j\subset {\cal F}$ denote the locus
$${\cal F}_j :=\left\{ ~Z\in {\cal F}~~|~~h^0(Z, \Omega^1_Z
\otimes \kappa^{-1/(n+1)})
\geq j~\right\} $$
where the space of candidate contact forms has dimension at least $j$.
 Since $\Omega^1_Z\otimes \kappa^{-1/(n+1)}$ is just the restriction
of $V^{\ast}\otimes L$ to the appropriate fiber of $\varpi$, and since $\varpi$
is
a flat morphism, it follows from the  semi-continuity theorem
\cite{hartsh} that each ${\cal F}_j$ is a Zariski-closed
subset of the quasi-projective variety
$\cal F$, and so, in particular, has only finitely many
components.

 Let $\tilde{\cal F}_j:= {\cal F}_j-{\cal F}_{j-1}$, $j\geq 1$. Since
$\tilde{\cal F}_j$ is a quasi-projective variety, it is a finite
union of  irreducible strata ${\cal F}_{jk}$,
each of which is a connected complex manifold.
On each stratum ${\cal F}_{jk}$, define a vector bundle ${ E}_{jk}$
as the zero-th direct image ${\cal O}(E_{jk})=\varpi^0_{\ast}{\cal O}
(V^{\ast}\otimes L)$
of the fiber-wise 1-forms with values in $L$.
Let $\cal L$ denote the line bundle
$\varpi^0_{\ast}{\cal O}(L^{n+1}\otimes \wedge^{2n+1}V)$
on $\cal F$. Then $\theta \mapsto \theta\wedge (d\theta)^{n}$ is
defines a canonical holomorphic section of the symmetric-product
bundle ${\cal L}\otimes \bigodot^{n+1}E_{jk}^{\ast }$; let
${\cal E}_{jk}$ denote the open subset in the total space of
$E_{jk}\to {\cal F}_{jk}$
where this homogeneous function is non-zero. Thus each
${\cal E}_{jk}$ is either a connected complex manifold  or is empty.
Each ${\cal E}_{jk}$ may now be viewed as the smooth parameter space of a
connected
family of Fano contact manifolds by taking the fiber over
 $\theta\in \Gamma (Z, \Omega^1(L))$ to be the pair $(Z,\ker \theta)$.
On the other hand, every  Fano contact manifold $(Z, D)$, where $Z$ is of the
fixed deformation type, appears in one of these families--- albeit many times.
Applying Lemma \ref{rig}, each of these families is of constant contact
type. Since we must  construct  $\cal F$
only for a finite number of degrees in order to account
for all Fano  deformation
types of the given dimension,
and since, for each $\cal F$ we only
have a finite number of contact families ${\cal E}_{jk}$,
 the result now  follows.
\end{proof}

\setcounter{main}{0}
\begin{main}[Finiteness Theorem]
 Up to  homothety,     there are
 only finitely many compact quaternionic-K\"ahler manifolds of
 positive scalar curvature in any given dimension $4n$.
\end{main}
\begin{proof}
By Proposition \ref{perfect}, two positive  quaternionic-K\"ahler
manifolds are homothetic iff their twistor spaces are
biholomorphic. Since the twistor space of any such manifold is
a Fano contact manifold, the result now follows immediately from
Theorem \ref{fanite}.\end{proof}

\section{Other Results}
We have seen in \S \ref{this} that the second homology of a
positive quaternionic-K\"ahler
manifold is far from arbitrary, and
may by itself contain enough information to determine the metric up to
isometry.
Recent calculations of Salamon show that the
higher homology groups are similarly constrained, in the
 following remarkable
manner:

\begin{thm}[Salamon]
 Let $(M^{4n},g)$ be a compact quaternionic-K\"ahler manifold
with positive scalar curvature. Then the
``odd''  Betti numbers $b_{2k+1}$ of $M$ vanish, and the ``even'' Betti numbers
$b_{2k}=b_{2(2n-k)}$ are subject to the linear constraint
$$   \sum_{k=0}^n  a_k  b_{2k}   =  0 ~,$$
where $a_k = \left\{\begin{array}{ll}1 + 2k + 2k ^2 - 4n/3 - 2kn + n^2 /3
&k < n\\(n^2  - n)/6 &k=n.
\end{array}\right. $
 \end{thm}
The proof of this result involves an intricate interplay between
the Kodaira vanishing theorem and the Penrose transform.
Details will appear elsewhere \cite{ls}.

\bigskip

\noindent {\bf Acknowledgements.}  The author would like to thank
Shigeru Mukai,
Jano\v{s} Koll\'ar, and Alan Nadel for their helpful explanations
of Fano theory, and Simon Salamon for many, many helpful conversations.


\begin{thebibliography}{99}

\bibitem{a2} D.V. Alekseevskii, {\em Compact Quaternion Spaces},
{\bf Funct.\ Anal.\ Appl.\ 2} (1968) 106-114.

\bibitem{a} D.V. Alekseevskii, {\em Classification of Quaternionic Spaces with
Transitive Solvable Groups of Motions}, {\bf Math.\ USSR--- Izv. 9} (1975)
297-339.

\bibitem{BE} T.N. Baily and M.G. Eastwood,
{\em Complex Paraconformal Manifolds---
Their Differential Geometry and Twistor Theory}, {\bf
 Forum Math.\ 3} (1991) 61-103.

\bibitem{BM} S. Bando and T. Mabuchi, {\em Uniqueness
of K\"ahler-Einstein Metrics Modulo Connected Group Actions}, {\bf
Algebraic Geometry, Sendai, 1985}, (T. Oda, ed), North-Holland 1985.

\bibitem{beber} L. B\'erard-Bergery, {\em Vari\'et\'es Quaternioniennes},
unpublished lecture notes, Espalion, 1979.

\bibitem{besse} A. Besse, {\bf Einstein Manifolds}, Springer-Verlag, 1987.


\bibitem{g} K. Galicki,{\em Generalization of the Momentum
Mapping Construction for Quaternionic K\"ahler Manifolds},
{\bf Comm.\ Math.\ Phys.\ 108} (1987) 117-138.

\bibitem{gl} K. Galicki and H.B. Lawson, {\em Quaternionic reduction and
quaternionic orbifolds,} {\bf Math.\ Ann.\  282} (1988) 1-21.

\bibitem{hartsh} R. Hartshorne, {\bf Algebraic Geometry}, Springer-Verlag,
1977.

\bibitem{H}  N.J. Hitchin, {\em K\"ahlerian Twistor Spaces},
{\bf  Proc.\ Lond.\ Math.\ Soc.\  43} (1981) 133-150.

\bibitem{kmm} J. Koll\'ar, Y. Miyaoka and S. Mori,
{\em Rational Connectedness and Boundedness of Fano manifolds},
{\bf J. Diff. Geom. 36} (1992) 765-779.

\bibitem{L1} C.R. LeBrun, {\em A Rigidity Theorem for Quaternionic-K\"ahler
  Manifolds}, {\bf Proc.\ Am.\ Math.\ Soc.\ 103} (1988) 1205-1208.

\bibitem{L0} C.R. LeBrun, {\em Quaternionic-K\"ahler  Manifolds
  and    Conformal   Geometry}, {\bf Math.\ Ann.\ 284} (1989) 353-376.

\bibitem{L2} C.R. LeBrun, {\em On Complete Quaternionic-K\"ahler Manifolds},
{\bf Duke Math.\ J. 63} (1991)
723-743.

\bibitem{L3} C.R. LeBrun, {\em On the Topology of Quaternionic Manifolds},
{\bf Twistor Newsletter 32} (1991) 6-7.

\bibitem{ls} C.R. LeBrun and S. Salamon, {\em
Strong Rigidity of Positive Quaternionic-K\"ahler Manifolds},
{\sl to appear}.

\bibitem{M}  S. Mori, {\em Hartshorne Conjecture and Extremal Ray},
{\bf Sugaku Expositions 0} (1988) 15-37.

\bibitem{nad} A. Nadel, {\em The Boundedness of Degree
of Fano Varieties with Picard Number One}, {\bf  J. Am.\ Math.\
Soc.\ 4} (1991) 681-692.

\bibitem{PS} Y.-S. Poon and S.M. Salamon,
{\em Eight-Dimensional Quaternionic-K\"ahler Manifolds with Positive Scalar
Curvature}, {\bf  J.\ Diff.\ Geometry }, {\em to appear}.

\bibitem{S} S.M. Salamon, {\em Quaternionic-K\"ahler Manifolds},
{\bf  Inv.\ Math.\ 67} (1982) 143-171.

\bibitem{Sbook} S.M. Salamon, {\bf Riemannian Geometry and Holonomy Groups},
Pitman Research Notes in Mathematics 201, Longman Scientific, 1989.

\bibitem{W} J.A. Wi\'sniewski, {\em On Fano Manifolds of Large Index},
{\bf Manu.\ Math.\ 70} (1991) 145-152.

\bibitem{wo} J.A. Wolf, {\em Complex Homogeneous Contact Manifolds and
Quaternionic Symmetric Spaces}, {\bf J. Math.\ Mech.\ 14} (1965) 1033-1047.


\end{thebibliography}
\end{document}